\documentclass{amsart}

\usepackage{amsmath,amssymb}
\usepackage{amsthm}
\usepackage{mathtools}

\DeclarePairedDelimiterX\innerp[2]{\langle}{\rangle}{#1, #2}
\DeclarePairedDelimiterX\ccint[2]{\lbrack}{\rbrack}{#1, #2}

\newtheorem{theorem}{Theorem}

\newtheorem{corollary}{Corollary}
\newtheorem{remark}{Remark}

\begin{document}

\bibliographystyle{amsplain}
\title{$L^p$ norm of truncated Riesz transform and an improved dimension-free $L^p$ estimate for maximal Riesz transform}

\author{Jinsong Liu}
\address{Jinsong Liu, HLM, Academy of Mathematics and Systems Science,
Chinese Academy of Sciences, Beijing, 100190, China \& School of Mathematical Sciences, University of Chinese Academy of Sciences, Beijing 100049,
China}
\email{liujsong@math.ac.cn}

\author{Petar Melentijevi\'c}
\address{Petar Melentijevi\'c, Matemati\v cki fakultet, University of Belgrade, Serbia}
\email{petar.melentijevic@matf.bg.ac.rs}

\author{Jian-Feng Zhu}
\address{Jian-Feng Zhu, School of Mathematical Sciences,
Huaqiao University,
Quanzhou 362021,
China.}
\email{flandy@hqu.edu.cn}

\date{}

\begin{abstract}
In this paper, we prove that the $L^p(\mathbb{R}^d)$ norm of the maximal truncated Riesz transform in terms of the $L^p(\mathbb{R}^d)$ norm of Riesz transform is dimension-free for any $2\leq p<\infty$, using integration by parts formula for radial Fourier multipliers. Moreover, we show that
$$\|R_j^*f\|_{L^p}\leq \left({2+\frac{1}{\sqrt{2}}}\right)^{\frac{2}{p}}\|R_jf\|_{L^p},\ \ \mbox{for}\ \  p\geq2,\ \  d\geq2.$$
As by products of our calculations, we infer the $L^p$ norm contractivity of the truncated Riesz transforms $R^t_j$ in terms of $R_j$, and their accurate $L^p$ norms. More precisely, we prove:
$$\|R^t_jf\|_{L^p}\leq\|R_jf\|_{L^p}$$
and
$$\|R^t_j\|_{L^p}=\|R_j\|_{L^p},$$
for all $1<p<+\infty,$ $j\in \{1,\dots,d\}$ and $t>0.$
\end{abstract}

\maketitle \pagestyle{myheadings} \markboth{Liu, Petar, and Zhu}{Dimension free estimate for the $L^p$ norm of $R_j^*$}


\section{Introduction}\label{sec-1}
Throughout the paper, we use $f*g$ to denote the convolution of $f$ and $g$, and use $\widehat{f}$ and $f^{\vee}$ to denote the {\it Fourier transform} and {\it inverse Fourier transform} of $f$, respectively. Moreover, we will use $\mathcal{M}f$ to denote the {\it centered Hardy-Littlewood maximal function}. Some special functions like {\it Bessel function} will be denoted by $J_{\alpha}(x)$, and the {\it Hypergeometric function} will be denoted by $ _pF_q$.

\subsection{Fourier and inverse Fourier transform}
Roughly speaking, a function is {\it Schwartz} if it is smooth and all of its derivatives decay faster than the reciprocal of any polynomial at infinity. For the precise definition we refer to \cite[Definition 2.2.1]{LG}.

Given a Schwartz function $f$ in $\mathbb{R}^d$, i.e. $f\in\mathcal{S}(\mathbb{R}^d)$, we define
$$\widehat{f}(\xi)=\int_{\mathbb{R}^d}f(x)e^{-2\pi i\langle x,\,\xi\rangle}\mathrm{d}x,$$
where $\langle x,\, \xi\rangle$
is the Euclidean inner product of $x$ and $\xi$.
We call $\widehat{f}$ the {\it Fourier transform} of $f$ (\cite[Definition 2.2.8]{LG}).

For all $x\in\mathbb{R}^d$, define
$$f^{\vee}(x)=\widehat{f}(-x).$$
The operator $f\mapsto f^{\vee}$ is called the {\it inverse Fourier transform} (\cite[Definition 2.2.13]{LG}).

\subsection{The maximal truncated Riesz transform}
To define the Riesz transform, we first introduce tempered distributions $\mathbf{W}_j$ on $\mathbb{R}^d$, for $1\leq j\leq d$, as follows: For $\varphi\in S(\mathbb{R}^d)$, let
$$\langle \mathbf{W}_j, \varphi\rangle=\frac{\Gamma(\frac{d+1}{2})}{\pi^{\frac{d+1}{2}}}\lim_{\varepsilon\to 0}\int_{|y|\geq\varepsilon}\frac{y_j}{|y|^{d+1}}\varphi(y)\mathrm{d}y.$$
The $j$th {\it Riesz transform} of $f$ is given by the convolution with the distribution $\mathbf{W}_j$, that is,
$$R_jf(x)=(f*\mathbf{W}_j)(x)=\frac{\Gamma(\frac{d+1}{2})}{\pi^{\frac{d+1}{2}}}\lim_{t\to 0^+}\int_{|x-y|>t}\frac{x_j-y_j}{|x-y|^{d+1}}f(y)\mathrm{d}y,$$
for all $f\in S(\mathbb{R}^d)$.

The above definition makes sense for any integrable function $f$ that has the property that for all $x$ there exist $C_x>0$, $\varepsilon_x>0$, and $\delta_x>0$ such that
for $y$ satisfying $|y-x|<\delta_x$, we have $|f(x)-f(y)|\leq C_x|x-y|^{\varepsilon_x}$. Moreover, one can show the following $L^p$ estimate:
$$\|R_jf\|_{L^p}\leq C\|f\|_{L^p},\ \ \ 1<p<\infty,$$
for some positive constant $C$ independent of $f$. Therefore, by the obvious density argument, $R_j$ can be defined on $L^p$ space (see \cite{Du}, \cite{MA-18}, and \cite{St}).

The {\it truncated Riesz transform}, denoted by $R_j^tf$, is given as follows:
\begin{equation}
\label{riesz-operator-Rj}R_j^tf(x)=\frac{\Gamma(\frac{d+1}{2})}{\pi^{\frac{d+1}{2}}}\int_{|x-y|>t}\frac{x_j-y_j}{|x-y|^{d+1}}f(y)\mathrm{d}y.
\end{equation}
It is easy to see that $R_jf(x)=\lim_{t\to 0^+}R_j^tf(x)$.

The {\it maximal truncated Riesz transform}, denoted by $R_j^*f$, is the maximal function of the singular integral (\ref{riesz-operator-Rj}). Namely, for $j=1, \ldots, d$, we set
$$R_j^*f(x)=\sup_{t>0}\left|R_j^tf(x)\right|.$$

\subsection{Norm estimate of the maximal truncated Riesz transform}
Clearly, $R_j^*f(x)\geq R_jf(x)$ for all $f\in \mathcal{S}(\mathbb{R}^d)$. On the other hand, Mateu and Verdera proved in \cite[Theorem 1]{MA-18} that, up to a multiplicative constant, also the reverse inequality holds in the $L^p(\mathbf{R}^d)$ norm. More precisely, they showed that there exists a constant $C_{p, d}$ depending on $p$ and the dimension $d$, such that
\begin{equation}
\label{norm-est-822}
\|R_j^*f\|_{L^p}\leq C_{p, d}\|R_jf\|_{L^p},
\end{equation}
for all $f\in L^p(\mathbb{R}^d)$ with $1<p<\infty$.
The above norm estimate (\ref{norm-est-822}) was later generalized to much broader classes of singular integral operators with even kernels \cite{MA-17} and with odd kernels \cite{MA-16}.

In 2022, Kucharski and Wr\'obel showed in \cite[Theorem 1.1]{KW-MA-22} that for the case $p=2$,  one may take a universal constant in (\ref{norm-est-822}), i.e. the $L^2$ norm of $R^*_j$ is dimension-free. More precisely, they proved the following result: For every $f\in L^2(\mathbb{R}^d)$, we have
$$\|R^*_jf\|_{L^2}\leq 2\cdot 10^8\|R_jf\|_{L^2}.$$
Very recently, they further improved the above estimate with Zienkiewicz, and obtain the dimension-free $L^p$ estimate for the vector of
the maximal Riesz transforms of odd order in terms of the Riesz transforms (see \cite[Theorem 2.3]{Wrobel-arxiv}) and then the same for all higher order Riesz transforms (see \cite[Theorem 2.3]{Wrobel-arxiv-2}).

\subsection{Two operators $M^t$ and $M^*$}
Let $M^t$, $t>0$, be defined by (see \cite[(3.1) and (3.5)]{KW-MA-22})
$$\widehat{M^tf}(\xi)=m(t|\xi|)\widehat{f}(\xi),\ \ \ f\in L^2(\mathbb{R}^d),\ \ \ \xi\in\mathbb{R}^d,$$
where
\begin{equation}
\label{function-mx}m(t|\xi|)=\frac{2^{\frac{d}{2}}\Gamma(\frac{d+1}{2})}{\sqrt{\pi}}\int_{2\pi t|\xi|}^{\infty}r^{-\frac{d}{2}}J_{\frac{d}{2}}(r)\mathrm{d}r,
\end{equation}
and set
$$M^*f(\xi)=\sup_{t>0}|M^tf(\xi)|.$$
Then for each $t>0$ and $f\in L^2(\mathbb{R}^d)$, the truncated Riesz transform factorizes as (\cite[(3.10)]{KW-MA-22})
\begin{equation}
\label{12-30-mt}R^t_jf=M^t(R_jf),\ \ \ j=1, \ldots, d.
\end{equation}
Moreover, the maximal operator $M^*$ is bounded on all $L^p(\mathbb{R}^d)$ spaces, $1<p<\infty$, and the optimal constant $C_p$ in the inequality
\begin{equation}
\label{conj-822}\|R^*_jf\|_{L^p}\leq C_p\|R_jf\|_{L^p}
\end{equation}
equals $\|M^*\|_{L^p\to L^p}$ (\cite[Corollary 3.3]{KW-MA-22}).

We will invoke (\ref{12-30-mt}) and (\ref{conj-822}) and the operators $M^t$ and $M^*$ in the rest of this paper. We remark here that in \cite[Theorem 2.3]{Wrobel-arxiv}, Kucharski, Wr\'obel and Zienkiewicz showed there exists a constant $C>0$, such that $C_p\leq C\cdot (p^*)^{2+1/2}$, where $p^*=\max\{p, (p-1)^{-1}\}$. 
\subsection{Motivation}
It would be interesting to investigate further the asymptotics of $C_p$ in (\ref{conj-822}) and more accurate order of its magnitude. Estimating maximal operators is often a difficult task and the more usual approach is using Littlewood-Paley decomposition and finding estimates of pieces whose Fourier transforms are well localized. Kucharski and Wrobel obtained their $L^2$ estimate using $M^*$ and following the Lemma 3 from \cite{Bourgain} with some subtle inequalities for Bessel functions. 
 
Let us say that the exact norms of maximal operators are even more delicate problems and are known only in a very limited number of cases. For example, weak type estimate for the centered Hardy-Littlewood maximal operator on $\mathbb{R}$ (see \cite{Ann-Math-03}) or $L^p$ inequality for uncentered Hardy-Littlewood maximal operator on $\mathbb{R}$ (see \cite{LG1}).

In this paper, we first give a sharp estimate of the function $m$ given by (\ref{function-mx}), and then, by using some identities for Bessel and hypergeometric functions, we find its inverse Fourier transform. This enables us to give the explicit formula of $M^t$. In fact, we split $M^tf$ into the sum of two convolution integrals $\varphi_1^t*f+\varphi_2^t*f$, where the kernel of $\varphi_1^t$ is supported away from the origin in $\{|x|\geq t\}$, while the kernel of $\varphi_2^t$ is supported near the origin in $\{|x|<t\}$. After precise calculations, we use the methods developed in \cite{AUSCHERCARO}. This approach gives better asymptotics for $C_p$ at infinity and considerably shorter proof. However, we were not been able to prove the same for $1<p<2.$

As an application of the obtained results, we show that $\|R_j^tf\|_{L^p}\leq \|R_jf\|_{L^p}$, for all $f\in L^p(\mathbb{R}^d)$ and $1<p<\infty$.
Moreover, using the same arguments developed in \cite{AUSCHERCARO} but for the concrete positive and radial kernel, we find a better dimension-free $L^p$ norm estimates for $R_j^*$ in terms of $R_j$.

Our main results are given in the next subsection. First we will give the accurate calculations for $m$ and its (inverse) Fourier transform, and then the $L^p$ estimates for truncated Riesz transforms and maximal operator $M^*.$

\subsection{Main results}
\begin{theorem}
\label{estimate-mx}
Let $m(t|\xi|)$ be the function given by $(\ref{function-mx})$. Then
$$1-\frac{2}{\pi}\int_0^{\pi}\frac{\sin t}{t}\mathrm{d}t\leq m\leq1,$$
and
$$\max_{0\leq t|\xi|<\infty}|m(t|\xi|)|=m(0)=1.$$
\end{theorem}

\begin{theorem}\label{Inverse-Fourier-m}
The inverse Fourier transform of $\xi \mapsto  m(t|\xi|)$ is as follows:
$$\varphi_t(x)=
\begin{cases}
\frac{\Gamma\left(\frac{1+d}{2}\right)^2t}{\pi^{1+\frac{d}{2}}|x|^{1+d}\Gamma\left(1+\frac{d}{2}\right)}\ _2F_1\left[\frac{1}{2}, \frac{1+d}{2}; \frac{2+d}{2}; \frac{t^2}{|x|^2}\right],\quad  \text{for} \quad t<|x|; \vspace{3mm}\\
\frac{\Gamma\left(\frac{1+d}{2}\right)^2t^{-d}}{\pi^{1+\frac{d}{2}}\Gamma\left(1+\frac{d}{2}\right)}\ _2F_1\left[\frac{1}{2}, \frac{1+d}{2}; \frac{2+d}{2}; \frac{|x|^2}{t^2}\right], \quad\quad \quad \text{for}\quad |x|<t.
\end{cases}
$$
Moreover, we have $\|\varphi_t\|_{L^1(\mathbb{R}^d)}=1$.
\end{theorem}
\begin{remark}
We remark here that in Appendix \ref{integral-37}, we have further proved that $m(t|\xi|)=O(|\xi|^{-(d+1)/2})$ as $|\xi|\to \infty$, see (\ref{9-16}). This shows $m\in L^2(\mathbb{R}^d)$, and thus, its Fourier transform and inverse Fourier transform are well defined in $L^2$ sense. In Theorem \ref{Inverse-Fourier-m}, for $|x|\neq t$, we calculate the inverse Fourier transform of $m$ in the sense of improper integral, and show it was $\varphi_t$. Then, by Plancherel's theorem, one has $\varphi_t\in L^2(\mathbb{R}^d)$.
On the other hand, since $\varphi_t(x)= O(|x|^{-d-1})$, as $|x|\to \infty$, one has $\varphi_t\in L^1(\mathbb{R}^d)$. In fact, we show more in Appendix \ref{L1-norm-phi} that $\|\varphi_t\|_{L^1(\mathbb{R}^d)}=1$. According to \cite[Theorem 9.13]{Rudin}, we have $m=\widehat{\varphi_t}$ and $M^tf=\varphi_t\ast f$.
\end{remark}

Using Theorem \ref{estimate-mx}, one can easily obtain the following corollary.

\begin{corollary}
\label{8-2-cor2}
Let $f\in L^2(\mathbb{R}^d)$. Then
$$\|R^t_jf\|_{L^2}^2\leq\|R_jf\|_{L^2}^2,\ \ \ j=1, \ldots, d.$$
\end{corollary}

It should be noted that Corollary \ref{8-2-cor2} is not clear from kernel definitions of the operators (as Professor Wr\'obel told us about this).
Moreover, as an application of Theorem \ref{estimate-mx} and Theorem \ref{Inverse-Fourier-m}, by using Young's inequality, one can easily find the following more general result.

\begin{theorem}
\label{HYtheorem}
Let $f\in L^p(\mathbb{R}^d),$ where $1<p<\infty$. Then
$$\|R^t_jf\|_{L^p}\leq\|R_jf\|_{L^p},\ \ \ j=1, \ldots, d,$$
and
$$\|R^t_j\|_{L^p}=\|R_j\|_{L^p}=H_p,$$
where
$$H_p=
\begin{cases}
\tan\frac{\pi}{2p},\quad  \text{for} \quad 1<p\leq2, \\
\cot\frac{\pi}{2p},\quad  \text{for}\quad 2\leq p<\infty,
\end{cases}
$$
was given by Iwaniec and Martin in \cite{IWANIECMARTIN}.
\end{theorem}

Finally, we give the dimension-free $L^p$ ($p\geq2$) norm estimate of $R^*_j$ in terms of $R_j$ as follows.

\begin{theorem}
\label{lp-estimate}
For every $f\in L^p(\mathbb{R}^d)$ with $p\geq2$, we have
\begin{equation}
\label{430-eq1}\|R_j^*f\|_{L^p}\leq \left(2+\frac{1}{\sqrt{2}}\right)^{\frac{2}{p}}\|R_jf\|_{L^p},\ \ \mbox{for}\ \  d\geq2.
\end{equation}
\end{theorem}


The proofs of Theorem \ref{lp-estimate}, Theorem \ref{estimate-mx}, Corollary \ref{8-2-cor2}, and Theorem \ref{HYtheorem} are given at Section \ref{sec-3}.

\section{Preliminaries}\label{sec-2}
In this section, we introduce two kind of special functions and recall two theorems which will be used in the proof of our main results. We start with the definition of the following Bessel function and hypergeometric function.

\subsection{Bessel functions}
Assume that $\alpha$ is a real number. The {\it Bessel function of first kind} $J_\alpha(x)$ is defined as follows (\cite[Page 200]{Landrews} or \cite[Page 15]{watson}):
\begin{equation}
\label{Bessel-function}J_\alpha(t)=\sum\limits_{k=0}^{\infty}\frac{(-1)^k}{\Gamma(k+\alpha+1)k!}\left(\frac{t}{2}\right)^{2k+\alpha}.
\end{equation}

It follows from \cite[(4.7.5)-(4.7.6)]{Landrews} that $J_\alpha(t)$ also has the following integral form:
\begin{align}\label{7-3-eq1}\nonumber
  J_\alpha(t) & =\frac{t^{\alpha}}{2^\alpha\Gamma(\alpha+\frac{1}{2})\sqrt{\pi}}\int_{-1}^1e^{its}(1-s^2)^{\alpha-\frac{1}{2}}\mathrm{d}s \\
   & =\frac{(t/2)^{\alpha}}{\sqrt{\pi}\Gamma(\alpha+\frac{1}{2})}\int_0^{\pi}\cos(t\cos\theta)\sin^{2\alpha}\theta\mathrm{d}\theta.
\end{align}

\subsection{Hypergeometric functions}
The hypergeometric function
$${}_pF_q[a_1, a_2, \ldots, a_p; b_1, b_2, \ldots, b_q; x]$$
is defined by the series (\cite[(2.1.2)]{Landrews}):
\begin{equation}
\label{Hypergeom-function}
{}_pF_q[a_1, a_2, \ldots, a_p; b_1, b_2, \ldots, b_q; x]=\sum\limits_{n=0}^{\infty}\frac{(a_1)_n\cdots(a_p)_n}{(b_1)_n\cdots(b_q)_n}\frac{x^n}{n!}
\end{equation}
for all $|x|<1$ and by continuation elsewhere.

Here $(q)_n$ is the Pochhammer symbol which is defined as follows
$$(q)_n=\left\{
\begin{array}
{r@{\ }l}
1, \ \  & \mbox{if}\ \ \ n=0;\\
\\
q(q+1)\cdots(q+n-1), \ \    & \mbox{if}\ \ \  n>0.
\end{array}\right.$$
Further we will use some formulas for hypergeometric and Bessel functions that can be found in \cite{HB-10}.

After making precise calculations, we will use the following general maximal estimate for radial multipliers:

\begin{theorem}\label{AC-11}$($\cite[Theorem 11 and Lemma 8]{AUSCHERCARO}$)$
Let $g$ be a radial kernel in $L^1(\mathbb{R}^d)$. Then, for $d\geq4$ and $p\geq2$,
$$\left\|\sup_{t>0}|g_t*f|\right\|_{L^p}\leq \left(C(2)+\sqrt{\frac{d-2}{d-3}}\right)^{2/p}\|g\|_{L^1}\|f\|_{L^p},$$
where $C(p)=p/(p-1).$
\end{theorem}

\section{Proofs of the main results}\label{sec-3}
\subsection{Proof of Theorem \ref{estimate-mx}}
For simplicity, let $x=t|\xi|$ and let $m(x)$ be the real-valued function defined by (\ref{function-mx}), i.e.
$$m(x)=\frac{2^{\frac{d}{2}}\Gamma(\frac{d+1}{2})}{\sqrt{\pi}}\int_{2\pi x}^{\infty}r^{-\frac{d}{2}}J_{\frac{d}{2}}(r)\mathrm{d}r,\ \ \ 0\leq x<\infty,$$
where $d\geq1$ is an integer.
It follows from \cite[10.22.43]{HB-10} that
$$m(0)=\frac{2^{\frac{d}{2}}\Gamma(\frac{d+1}{2})}{\sqrt{\pi}}\int_{0}^{\infty}r^{-\frac{d}{2}}J_{\frac{d}{2}}(r)\mathrm{d}r=1.$$
Thus, by using (\ref{7-3-eq1}), we have
\begin{align*}
  m(x) & =1-\frac{2^{\frac{d}{2}}\Gamma(\frac{d+1}{2})}{\sqrt{\pi}}\int_{0}^{2\pi x}r^{-\frac{d}{2}}J_{\frac{d}{2}}(r)\mathrm{d}r \\
   & =1-\frac{1}{\pi}\int_0^{\pi}\sin^{d}\theta\bigg(\int_0^{2\pi x}\cos(r\cos\theta)\mathrm{d}r\bigg)\mathrm{d}\theta\\
   &=1-\frac{2}{\pi}\int_0^{\frac{\pi}{2}}\frac{\sin^d\theta\sin(2\pi x\cos\theta)}{\cos \theta}\mathrm{d}\theta.
\end{align*}
Let
$$I(x)=\frac{2}{\pi}\int_0^{\frac{\pi}{2}}\frac{\sin^d\theta\sin(2\pi x\cos\theta)}{\cos \theta}\mathrm{d}\theta=1-m(x).$$
Then, by letting $t=2\pi x\cos \theta$, we have
\begin{align}\label{8-13-eq1}\nonumber
  I(x) & =\frac{2}{\pi}\int_0^{\frac{\pi}{2}}\sin^{d-1}\theta\frac{\sin(2\pi x\cos\theta)}{2\pi x\cos\theta}2\pi x\sin\theta\mathrm{d}\theta \\
    & =\frac{2}{\pi}\int_0^{2\pi x}\left(1-\left(\frac{t}{2\pi x}\right)^2\right)^{\frac{d-1}{2}}\frac{\sin t}{t}\mathrm{d}t.
\end{align}

To estimate $I(x)$, since $0\leq1-(t/(2\pi x))^2\leq1$, one may focus on the sign of $\sin t$.
First, we show $I(x)\geq0$ as follows:
If $0\leq x\leq1/2$, then $I(x)\geq0$ holds trivial.
If $1/2<x<1$, then $0<2\pi x-\pi<\pi$, and we have
$$I(x)=\frac{2}{\pi}\int_0^{\pi}\left(1-\left(\frac{t}{2\pi x}\right)^2\right)^{\frac{d-1}{2}}\frac{\sin t}{t}\mathrm{d}t-\frac{2}{\pi}\int_0^{2\pi x-\pi}\left(1-\left(\frac{\pi+t}{2\pi x}\right)^2\right)^{\frac{d-1}{2}}\frac{\sin t}{\pi+t}\mathrm{d}t\geq0.$$
Now, let $k\geq1$ be an integer, and suppose $k\leq{x}<k+1$.
Then
\begin{align}\label{412-eq-2}\nonumber
  I(x) & =\frac{2}{\pi}\sum_{j=0}^{2k-1}\int_{j\pi}^{(j+1)\pi} \bigg(1-\frac{t^2}{4\pi^2x^2}\bigg)^{\frac{d-1}{2}}\frac{\sin t}{t}\mathrm{d}t+\frac{2}{\pi}\int_{2k\pi}^{2\pi x} \bigg(1-\frac{t^2}{4\pi^2x^2}\bigg)^{\frac{d-1}{2}}\frac{\sin t}{t}\mathrm{d}t \\
   & =\frac{2}{\pi}\sum_{j=0}^{2k-1}(-1)^j a_j(x)+\frac{2}{\pi}\int_{2k\pi}^{2\pi x} \bigg(1-\frac{t^2}{4\pi^2x^2}\bigg)^{\frac{d-1}{2}}\frac{\sin t}{t}\mathrm{d}t,
\end{align}
where
$$a_j(x)=\int_{0}^{\pi} \bigg(1-\frac{(t+j\pi)^2}{4\pi^2x^2}\bigg)^{\frac{d-1}{2}}\frac{\sin t}{t+j\pi}\mathrm{d}t.$$
Notice that $t+j\pi\leq 2k\pi\leq2\pi x$ holds for all $0\leq j\leq 2k-1$ and all $0\leq t\leq \pi$. We see that $a_j(x)$ is non-negative and evidently decreasing in $j$. Therefore,
$$\sum_{j=0}^{2k-1}(-1)^j a_j(x)=(a_0(x)-a_1(x))+(a_2(x)-a_3(x))+\dots+(a_{2k-2}(x)-a_{2k-1}(x))\geq 0.$$

To show the last integration of (\ref{412-eq-2}) is non-negative, we remark that if $k\leq x\leq k+1/2$, then the integral function is evidently non-negative; for another case $k+1/2<x<k+1$, we have
\begin{align*}
  \int_{2k\pi}^{2\pi x} \bigg(1-\frac{t^2}{4\pi^2x^2}\bigg)^{\frac{d-1}{2}}\frac{\sin t}{t}\mathrm{d}t & = \left(\int_{2k\pi}^{2k\pi +\pi}+ \int_{2k\pi+\pi}^{2\pi x}\right) \bigg(1-\frac{t^2}{4\pi^2x^2}\bigg)^{\frac{d-1}{2}}\frac{\sin t}{t}\mathrm{d}t \\
   &\geq a_{2k}(x)-a_{2k+1}(x)\geq 0.
\end{align*}

Next, we are going to find an upper bound for $I(x)$ as follows:
If $0\leq x<1/2$, then
$$I(x)=\frac{2}{\pi}\int_{0}^{2\pi x} \bigg(1-\frac{t^2}{4\pi^2x^2}\bigg)^{\frac{d-1}{2}}\frac{\sin t}{t}\mathrm{d}t \leq\frac{2}{\pi}\int_{0}^{\pi} \frac{\sin t}{t}\mathrm{d}t.$$
If $1/2\leq x<1$, then
$$I(x)=\frac{2}{\pi}\left(\int_{0}^{\pi}+\int_{\pi}^{2\pi x}\right) \bigg(1-\frac{t^2}{4\pi^2x^2}\bigg)^{\frac{d-1}{2}}\frac{\sin t}{t}\mathrm{d}t \leq\frac{2}{\pi}\int_{0}^{\pi} \frac{\sin t}{t}\mathrm{d}t.$$
Now, fix $k\geq1$ and suppose $k\leq x<k+1/2$. One has
$$\int_{2k\pi}^{2\pi x} \bigg(1-\frac{t^2}{4\pi^2x^2}\bigg)^{\frac{d-1}{2}}\frac{\sin t}{t}\mathrm{d}t \leq \int_{2k\pi}^{2k\pi +\pi} \bigg(1-\frac{t^2}{4\pi^2x^2}\bigg)^{\frac{d-1}{2}}\frac{\sin t}{t}\mathrm{d}t=a_{2k}(x).$$
Similarly, if $k+1/2\leq x<k+1$, then again, one has
\begin{align*}
  \int_{2k\pi}^{2\pi x} \bigg(1-\frac{t^2}{4\pi^2x^2}\bigg)^{\frac{d-1}{2}}\frac{\sin t}{t}\mathrm{d}t& =a_{2k}(x)+\int_{(2k+1)\pi}^{2\pi x} \bigg(1-\frac{t^2}{4\pi^2x^2}\bigg)^{\frac{d-1}{2}}\frac{\sin t}{t}\mathrm{d}t \\
   &\leq a_{2k}(x).
\end{align*}
Moreover, for $k\geq 1$ and $k\leq x<k+1$, one has
$$\sum_{j=0}^{2k-1}(-1)^j a_j(x)=a_0(x)-(a_1(x)-a_2(x))-\dots-(a_{2k-3}(x)-a_{2k-2}(x))-a_{2k-1}(x).$$
We see from (\ref{412-eq-2}) that
$$I(x)\leq\frac{2}{\pi} \left(\sum_{j=0}^{2k-1}(-1)^j a_j(x)+a_{2k}(x)\right)\leq\frac{2}{\pi}a_0(x),\ \ \ 1\leq x<\infty.$$
Combining all the cases $0\leq x<1/2$, $1/2\leq x<1$, and $1\leq x<\infty$, we have showed that
$$I(x)\leq\frac{2}{\pi}\int_{0}^{\pi} \frac{\sin t}{t}\mathrm{d}t,\ \ \ \mbox{for all}\ \ \ 0\leq x<\infty.$$

Based on the above discussions, we have
$$-0.17898\approx 1-\frac{2}{\pi}\int_0^{\pi}\frac{\sin t}{t}\mathrm{d}t\leq m(x)=1-I(x)\leq1,$$
which completes the proof.
\qed

\begin{remark}
\label{rem-822}
Let us explain some facts which lie in the above proof.
Let
$$\phi(t)=\frac{2}{\pi}\left(1-\left(\frac{t}{2\pi x}\right)^2\right)^{\frac{d-1}{2}}\frac{\sin t}{t},$$
where $0\leq t\leq2\pi x$ and $d\geq1$. Then
$$I(x)=\int_0^{2\pi x}\phi(t)\mathrm{d}t.$$

It is easy to see that
\begin{equation}
\label{717-1}\phi(t)\geq0,\ \ \ \mbox{if} \ \ t\in [2k\pi, (2k+1)\pi]\cap[0, 2\pi x],
\end{equation}
while
\begin{equation}
\label{717-2}\phi(t)\leq0,\ \ \ \mbox{if} \ \ t\in [(2k+1)\pi, (2k+2)\pi]\cap[0, 2\pi x],
\end{equation}
where $k$ are non-negative integers.

Moreover, for any $0\leq t\leq\pi+t\leq 2\pi x$, the following inequality holds
\begin{equation}
\label{7-8-eq2}|\phi(\pi+t)|\leq|\phi(t)|.
\end{equation}
The above facts (\ref{717-1}), (\ref{717-2}), and (\ref{7-8-eq2}) are in fact repeated many times in the proof of Theorem \ref{estimate-mx}.
\end{remark}

\subsection{Proof of Theorem \ref{Inverse-Fourier-m}}
By using the (power series) definition of Bessel functions (\eqref{Bessel-function}) and Hypergeometric functions (\eqref{Hypergeom-function}), we calculate the function $m(t|\xi|)$ as follows:
\begin{align}\label{12-29-mx}
  m(t|\xi|) & =1-\frac{2^{\frac{d}{2}}\Gamma(\frac{d+1}{2})}{\sqrt{\pi}}\int_{0}^{2\pi t|\xi|}r^{-\frac{d}{2}}J_{\frac{d}{2}}(r)\mathrm{d}r \\\nonumber
   & =1-\frac{\Gamma\left(\frac{1+d}{2}\right)\pi t|\xi|}{\Gamma\left(\frac{3}{2}\right)\Gamma\left(\frac{2+d}{2}\right)}\ _1F_2\left[\frac{1}{2}; \frac{3}{2}, 1+\frac{d}{2}; -\pi^2(t|\xi|)^2\right].
\end{align}
Let $\varphi_t=m^{\vee}(t|\xi|)$ be the inverse Fourier transform of $\xi \mapsto m(t|\xi|)$. Then
$$\varphi_t(x)=\int_{\mathbb{R}^d}m(t|\xi|)e^{2\pi i \langle x,\, \xi\rangle}\mathrm{d}\xi,$$
where $x, \xi\in\mathbb{R}^d$.

Since $m(t|\xi|)$ is a radial function, for $d\geq2$, we see from \cite[Appendix B.5]{LG} that $\varphi_t(x)$ can be rewritten as follows:
\begin{align}\label{appendix}\nonumber
  \varphi_t(x)=\int_{\mathbb{R}^d}m(t|\xi|)e^{2\pi i \langle x,\, \xi\rangle}\mathrm{d}\xi&=\int_{0}^{\infty}m(tr)r^{d-1}\bigg(\int_{\mathbb{S}^{d-1}}e^{2\pi i \langle x,\, r\theta\rangle}\mathrm{d}\theta\bigg)\mathrm{d}r \\
      &=\frac{2\pi}{|x|^{\frac{d-2}{2}}}\int_{0}^{\infty}m(tr)J_{\frac{d-2}{2}}(2\pi r |x|)r^{\frac{d}{2}}\mathrm{d}r.
\end{align}
In order to calculate the above integral, we need the following three formulas (see \cite[(4.6.1)]{Landrews}):
\begin{equation}
\label{1-08-eq1}r^\alpha J_{\alpha-1}(r)=\frac{\mathrm{d}}{\mathrm{d}r}r^\alpha J_{\alpha}(r),
\end{equation}
and (see \cite[10.22.56]{HB-10}):
\begin{equation}
\label{1-08-eq2}\int_0^{\infty}J_{\frac{d}{2}}(\rho)J_{\frac{d}{2}}\left(\frac{t\rho}{|x|}\right)\mathrm{d}\rho=\frac{\Gamma(\frac{1+d}{2})t^{\frac{d}{2}}{}_2F_1\left[\frac{1}{2}, \frac{1+d}{2}; \frac{2+d}{2}; \frac{t^2}{|x|^2}\right]}{\Gamma(1+\frac{d}{2})\sqrt{\pi}|x|^{\frac{d}{2}}}, \ t<|x|,
\end{equation}
and (by changing integral variable: $u={t\rho}/{|x|}$)
\begin{equation}
\label{4-26-eq1}\int_0^{\infty}J_{\frac{d}{2}}(\rho)J_{\frac{d}{2}}\left(\frac{t\rho}{|x|}\right)\mathrm{d}\rho=\frac{\Gamma(\frac{1+d}{2})t^{-\frac{d}{2}-1}{}_2F_1\left[\frac{1}{2}, \frac{1+d}{2}; \frac{2+d}{2}; \frac{|x|^2}{t^2}\right]}{\sqrt{\pi}\Gamma(1+\frac{d}{2})|x|^{-\frac{d}{2}-1}},\ t>|x|.
\end{equation}

Now, integral by parts in (\ref{appendix}), and using (\ref{12-29-mx}), (\ref{1-08-eq1}), (\ref{1-08-eq2}), one has
\begin{equation}
\label{1-9-eq1}\varphi_t(x)=\frac{\Gamma\left(\frac{1+d}{2}\right)^2t}{\pi^{1+\frac{d}{2}}|x|^{1+d}\Gamma\left(1+\frac{d}{2}\right)}\ _2F_1\left[\frac{1}{2}, \frac{1+d}{2}; \frac{2+d}{2}; \frac{t^2}{|x|^2}\right],\ \mbox{if}\ t<|x|; 
\end{equation}
and similarly, using (\ref{4-26-eq1}), one has
\begin{equation}
\label{4-26-Im-2} \varphi_t(x)=\frac{\Gamma\left(\frac{1+d}{2}\right)^2t^{-d}}{\pi^{1+\frac{d}{2}}\Gamma\left(1+\frac{d}{2}\right)}\ _2F_1\left[\frac{1}{2}, \frac{1+d}{2}; \frac{2+d}{2}; \frac{|x|^2}{t^2}\right],\  \mbox{if}\  t>|x|. 
\end{equation}
For more details of the above calculations, we refer to Appendix. Let us say that for the case $d=1$ there holds the same formula, the difference is that in the formula \eqref{appendix} we do not use the integral over the sphere but the fact that $m$ is even. 

Next, we calculate the $L^1$ norm of $\varphi_t$ as follows: Since $\varphi_t$ is a positive function, the condition $m(0)=1$ easily implies $\|\varphi_t\|_{L^1}=1$, since
$$m(0)=\widehat{\varphi}_t(0)=\int_{\mathbb{R}^d}\varphi_t(x)\mathrm{d}x=\|\varphi_t\|_{L^1}.$$
For the reader's convenience, we also give an another computational proof of $\|\varphi_t\|_{L^1}=1$ at Appendix.
\qed
\vspace{5mm}

The following proof was due to Professor Wr\'obel.
\subsection{Proof of Corollary \ref{8-2-cor2}}
Since $R^t_jf=M^t(R_jf)$ and  $\widehat{M^t(R_jf)}(\xi)=m(t|\xi|)\widehat{R_jf}(\xi)$, according to Theorem \ref{estimate-mx} and Plancherel's theorem, we have
\begin{align*}
  \|R^t_jf\|_{L^2}^2 &= \int_{\mathbb{R}^d}|m(t|x|)|^2|\widehat{R_jf(x)}|^2\mathrm{d}x \\
   &\leq\int_{\mathbb{R}^d}|\widehat{R_jf(x)}|^2\mathrm{d}x\\
   &=\|R_jf\|_{L^2}^2,
\end{align*}
which completes the proof. \qed

\subsection{Proof of Theorem \ref{HYtheorem}} According to the Young's inequality, we see from (\ref{12-30-mt}) and Theorem \ref{Inverse-Fourier-m} that:
$$\|R^t_jf\|_{L^p}=\|M^t(R_jf)\|_{L^p}=\|\varphi_t\ast R_jf\|_{L^p}\leq \|\varphi_t\|_{L^1}\|R_jf\|_{L^p}=\|R_jf\|_{L^p},$$
where $f\in L^p(\mathbb{R}^d)\cap L^2(\mathbb{R}^d)$ and $1<p<\infty$.

Let $f_n\in L^p(\mathbb{R}^d)\cap L^2(\mathbb{R}^d)$ such that $\lim_{n\to \infty}f_n=f\in L^p(\mathbb{R}^d)$. Then
\begin{align*}
  \|R_j^tf\|_{L^p} & \leq \|R_j^tf_n\|_{L^p}+\|R_j^t(f-f_n)\|_{L^p} \\
  &\leq\|R_jf_n\|_{L^p}+\|R_j^*(f-f_n)\|_{L^p}\\
   &\leq \|R_jf\|_{L^p}+\|R_j(f-f_n)\|_{L^p}+C_{p, d}\|R_j(f-f_n)\|_{L^p},
\end{align*}
where $C_{p, d}$ is a constant depending on $p$ and $d$ (see (\ref{norm-est-822})).
Letting $n\to \infty$, and noticing that $\|R_j(f-f_n)\|_{L^p}\leq C\|f-f_n\|_{L^p}$, we have
$$\|R_j^tf\|_{L^p}\leq\|R_jf\|_{L^p},\ \ \ 1<p<\infty.$$

It follows from \cite[Theorem 1.1]{IWANIECMARTIN} that $\|R_jf\|_{L^p}\leq H_p \|f\|_{L^p}$,
where
$$H_p=
\begin{cases}
\tan\frac{\pi}{2p},\quad  \text{for} \quad 1<p\leq2, \\
\cot\frac{\pi}{2p},\quad  \text{for}\quad 2\leq p<\infty.
\end{cases}
$$
Therefore, one can easily get $\|R^t_jf\|_{L^p}\leq H_p\|f\|_{L^p}$, and thus,
$$\|R^t_j\|_{L^p}\leq H_p$$
for all $t>0$.

On the other hand, recall that
$$R^t_jf(x)=\frac{\Gamma(\frac{d+1}{2})}{\pi^{\frac{d+1}{2}}}\int_{|x-y|>t}\frac{x_j-y_j}{|x-y|^{d+1}}f(y)\mathrm{d}y.$$
By substitution $x=tu$ and $y=tz$, we get
$$R^t_jf(tu)=\frac{\Gamma(\frac{d+1}{2})}{\pi^{\frac{d+1}{2}}}\lim_{t\to 0^+}\int_{|u-z|>1}\frac{u_j-z_j}{|u-z|^{d+1}}f(tz)\mathrm{d}z=R^1_jf_t(u),$$
where $f_t(u):=f(tu).$
Hence, the norm of $R^t_j$ does not depend on $t$, and therefore, by $\lim_{t \rightarrow 0^+}R^{t}_jf=R_jf$ and Fatou's lemma we have:
$$\|R^{t}_j\|_{L^p}\geq \left\|\lim_{t \rightarrow 0^+}R^{t}_j\right\|_{L^p}=\|R_j\|_{L^p}=H_p.$$
Based on the above discussions, we have
$$\|R^t_j\|_{L^p}=H_p$$
and the proof is completed. \qed

Let us say that this method of getting the estimate from below was first exploited by \cite{LAENG}.

\subsection{Proof of Theorem \ref{lp-estimate}}
Recall that the operator $M^t$, $t>0$, is defined by
$$\widehat{M^tf}(\xi)=m(t|\xi|)\widehat{f}(\xi), \ \ \ f\in \mathcal{S}(\mathbb{R}^d).$$
Therefore,
$$ M^tf(x)=\left(m(t|\xi|)\widehat{f}(\xi)\right)^{\vee}(x)=\int_{\mathbb{R}^d}m(t|\xi|)\widehat{f}(\xi)e^{2\pi i\langle x, \xi\rangle}\mathrm{d}\xi,$$
where
 $$\widehat{f}(\xi)=\int_{\mathbb{R}^d}f(y)e^{-2\pi i \langle y,\, \xi\rangle}\mathrm{d}y$$
is the Fourier transform of $f$.

Then
\begin{align*}
   M^tf(x)& =\int_{\mathbb{R}^d}\int_{\mathbb{R}^d}f(y)m(t|\xi|)e^{-2\pi i \langle y-x,\, \xi \rangle}\mathrm{d}\xi\mathrm{d}y \\
   & =\int_{\mathbb{R}^d}f(y)\bigg(\int_{\mathbb{R}^d}m(t|\xi|)e^{2\pi i \langle x-y,\, \xi\rangle}\mathrm{d}\xi\bigg)\mathrm{d}y\\
   &:=(\varphi_t\ast f)(x),
\end{align*}
where
$$\varphi_t(x)=\int_{\mathbb{R}^d}m(t|\xi|)e^{2\pi i \langle x,\, \xi\rangle}\mathrm{d}\xi$$
is the inverse Fourier transform of $m(t|\xi|)$ and was given in Theorem \ref{Inverse-Fourier-m}. 

Using these facts, one can find the concrete formula of $M^t$ as follows:
\begin{align*}
  M^tf(x) &=\frac{\Gamma\left(\frac{1+d}{2}\right)^2}{\pi^{1+\frac{d}{2}}\Gamma\left(1+\frac{d}{2}\right)}\int_{|y-x|>t}\frac{t}{|y-x|^{1+d}}\ _2F_1\left[\frac{1}{2}, \frac{1+d}{2}; \frac{2+d}{2}; \frac{t^2}{|y-x|^2}\right]f(y)\mathrm{d}y \\
   &\ \ \ +\frac{\Gamma\left(\frac{1+d}{2}\right)^2t^{-d}}{\pi^{1+\frac{d}{2}}\Gamma\left(1+\frac{d}{2}\right)}\int_{|y-x|<t}\ _2F_1\left[\frac{1}{2}, \frac{1+d}{2}; \frac{2+d}{2}; \frac{|y-x|^2}{t^2}\right]f(y)\mathrm{d}y.
\end{align*}
For simplicity, if we define $\varphi_1^{t}=\varphi_t(x)\chi_{|x|>1}(x)$ and $\varphi_2^t=\varphi_t(x)\chi_{|x|<1}(x)$, then
$$M^tf=\varphi_1^t*f+\varphi_2^t*f.$$

Next, we are going to estimate the $L^p$ norm of maximal operator $M^*f=\sup_{t>0}|M^tf|$.
Let
\begin{equation}
\label{function-phi}\varphi(x)=
\begin{cases}
\frac{\Gamma\left(\frac{1+d}{2}\right)^2}{\pi^{1+\frac{d}{2}}\Gamma\left(1+\frac{d}{2}\right)}\frac{1}{|x|^{d+1}}\ _2F_1\left[\frac{1}{2}, \frac{1+d}{2}; \frac{2+d}{2}; \frac{1}{|x|^2}\right],\quad  \text{for} \quad |x|>1, \vspace{3mm}\\
\frac{\Gamma\left(\frac{1+d}{2}\right)^2}{\pi^{1+\frac{d}{2}}\Gamma\left(1+\frac{d}{2}\right)}\ _2F_1\left[\frac{1}{2}, \frac{1+d}{2}; \frac{2+d}{2}; |x|^2\right],\quad\quad \quad \text{for}\quad |x|<1.
\end{cases}
\end{equation}
Then $\varphi_t(x):=\varphi(x/t)/t^d$, and $\|\varphi\|_{L^1}=\|\varphi_t\|_{L^1}=1$, since both $\varphi_t$ and $\varphi$ are positive and radial (We refer to Appendix for more detailed calculations on this). 

It is easy to see  that $\varphi$ is continuous and decreasing, when $|x|>1$; and $\varphi$ is increasing, when $|x|<1$. This enables us to use the usual approach of estimating maximal operator $M^*$ via least decreasing radial majorant and Hardy-Littlewood maximal operator. Instead, we use the Theorem A (\cite{AUSCHERCARO}). Using it directly as it was stated, one can easily obtain the following estimate: For $d\geq4$ and $p\geq2$,
\begin{align*}
\|M^*f\|_{L^p}&=\left\|\sup_{t>0}|M^tf|\right\|_{L^p}=\left\|\sup_{t>0}|\varphi_t*f|\right\|_{L^p}\\
&\leq (2+\sqrt{2})^{\frac{2}{p}}\|\varphi\|_{L^1}\|f\|_{L^p}=(2+\sqrt{2})^{\frac{2}{p}}|f\|_{L^p}.
\end{align*}

However, our main purpose is to find estimate for all $d\geq1$. Since we have calculated the kernel $\varphi_t$, following the arguments from \cite{AUSCHERCARO}, we obtain a new and better estimate for $\|M^*f\|_{L^p}$ when $d\geq2$. We start with the following identity, which follows from integration by parts:
$$\varphi_t\ast f=\frac{1}{t}\int_{0}^{t}\varphi_s\ast f \mathrm{d}s+\frac{1}{t}\int_{0}^{t}s\frac{d}{ds}\varphi_s\ast f \mathrm{d}s.$$
According to \cite[Proposition 9]{AUSCHERCARO}, the first term is majorized by $2\|f\|_{L^2}$, i.e.
$$\left\|\sup_{r>0}\left|\frac{1}{t}\int_{0}^{t}\varphi_s\ast f \mathrm{d}s\right|\,\right\|_{L^2}\leq 2\|\varphi\|_{L^1}\|f\|_{L^2}=2\|f\|_{L^2},$$
while the second term is bounded by
$$\sup_{t>0}\bigg|\frac{1}{t}\int_{0}^{t}s\frac{d}{ds}\varphi_s\ast f \mathrm{d}s\bigg|\leq \bigg(\int_{0}^{+\infty}\big|s\frac{d}{ds}\varphi_s\ast f\big|^2 \frac{\mathrm{d}s}{s}\bigg)^{\frac{1}{2}}:=\widetilde{g}(f).$$
Using Plancherel's theorem, the $L^2$ norm of the last Littlewood-Paley function can be estimated as
$$\left\|\widetilde{g}(f)\right\|_{L^2} \leq \sup_{\xi \in \mathbb{R}^d}\bigg(\int_{0}^{+\infty}\big|s\frac{d}{ds}\hat\varphi_s(\xi)\big|^2 \frac{\mathrm{d}s}{s}\bigg)^{\frac{1}{2}}\|f\|_{L^2}.$$
Since $\hat\varphi_s(\xi)=m(s|\xi|)$, we see from (\ref{12-29-mx}) and \cite[(10.22.57)]{HB-10} that
\begin{align*}
\int_{0}^{+\infty}\big|s\frac{d}{ds}\hat\varphi_s(\xi)\big|^2 \frac{\mathrm{d}s}{s}&=\frac{2^d\Gamma(\frac{d+1}{2})^2}{\pi}\int_{0}^{+\infty} \left((2\pi s|\xi|)^{1-\frac{d}{2}}J_{\frac{d}{2}}(2\pi s|\xi|)\right)^2\frac{\mathrm{d}s}{s}\\
&=\frac{2^d\Gamma(\frac{d+1}{2})^2}{\pi}\int_{0}^{+\infty} s^{1-d}J_{\frac{d}{2}}(s)^2\mathrm{d}s\\
&=\frac{2}{\pi}\frac{\Gamma(\frac{d+1}{2})^2}{(d-1)\Gamma(\frac{d}{2})^2}.
\end{align*}
Let us prove that the last expression is $\leq1/2$, for all $d\geq 2$, or what is the same, $$a_d=\frac{\Gamma(\frac{1}{2})\Gamma(\frac{d}{2})}{\Gamma(\frac{d+1}{2})}\sqrt{d-1}\geq 2, \ \ \ \mbox{for all}\ \ \ d\geq2.$$
But this easily follows from the following facts: ${a_{d+2}}/{a_d}=\sqrt{{d^2}/(d^2-1)}>1$ and $a_2=2, a_3={\pi}/{\sqrt{2}}>2$.
Then
$$\left\|\sup_{t>0}\left|\frac{1}{t}\int_{0}^{t}s\frac{d}{ds}\varphi_s\ast f \mathrm{d}s\right|\,\right\|_{L^2}\leq \frac{1}{\sqrt{2}}\|f\|_{L^2}.$$

Therefore, for $d\geq 2$ and $p\geq 2,$ using interpolation, we conclude
$$\|M^*f\|_{L^p} \leq \bigg(2+\frac{1}{\sqrt{2}}\bigg)^{\frac{2}{p}}\|f\|_{L^p},$$
i.e.
$$\|R_j^*f\|_{L^p} \leq \bigg(2+\frac{1}{\sqrt{2}}\bigg)^{\frac{2}{p}}\|R_jf\|_{L^p}.$$
Here, notice that for $p=\infty$, we use the fact $\|\varphi_t\|_{L^1}=\|\varphi\|_{L^1}=1$.

For the case of $d=1$, the maximal Riesz transform is just the Hilbert transform $H^*$ (see \cite[Definition 5.1.10]{LG}).
We then use the usual approach which was given in \cite{LG}. From the inequality \cite[(5.1.35)]{LG}
$$|H^*f(x)|\leq \frac{1}{\pi}\|\Psi\|_{L^1}\mathcal{M}(f)(x)+\mathcal{M}(H(f))(x),$$
where (\cite[(5.1.21)]{LG})
$$\Psi(x)=
\begin{cases}
\frac{1}{x(1+x^2)},\quad  \text{when} \quad |x|>1, \\
\frac{1}{2},\quad \quad \ \ \text{when}\quad |x|<1,
\end{cases}
$$
we get
$$\|H^*f\|_{L^2}\leq \frac{1}{\pi}\|\Psi\|_{L^1}\|\mathcal{M}f\|_{L^2}+\|\mathcal{M}(H(f))\|_{L^2}.$$
By using explicit constants for $L^1$ norm of $\Psi$, the best constant for uncentered Hardy-Littlewood maximal operator, and unitarity of Hilbert transform, we have
$$\|H^*f\|_{L^2}\leq (1+\sqrt{2})\bigg(1+\frac{1+\log 2}{\pi}\bigg)\|Hf\|_{L^2}<\frac{15}{4}\|Hf\|_{L^2}.$$
For $L^{\infty}$ estimate we use that $\|\varphi_t\|_{L^1}=1$ and representation of $M^t$ as convolution operator. Again, by interpolation and the $L^{\infty}$ estimate of $H^*$, we finally get:
$$\|H^*f\|_{L^p}<\bigg(\frac{15}{4}\bigg)^{\frac{2}{p}}\|Hf\|_{L^p},$$
for $p\geq 2.$
\qed

\begin{remark}
For the case of $d=1$, the maximal Riesz transform is just the maximal Hilbert transform $H^*$ (see \cite[Definition 5.1.10]{LG}).
We then use the usual approach which was given in \cite{LG}. From the inequality \cite[(5.1.35)]{LG}
$$|H^*f(x)|\leq \frac{1}{\pi}\|\Psi\|_{L^1}\mathcal{M}(f)(x)+\mathcal{M}(H(f))(x),$$
where $H$ is the Hilbert transform and (see \cite[(5.1.21)]{LG}, but our majorant here is a slightly smaller one)
$$\Psi(x)=
\begin{cases}
\frac{1}{x(1+x^2)},\quad  \text{when} \quad |x|>1, \\
\frac{1}{2},\quad \quad \ \ \text{when}\quad |x|<1,
\end{cases}
$$
we get
$$\|H^*f\|_{L^2}\leq \frac{1}{\pi}\|\Psi\|_{L^1}\|\mathcal{M}f\|_{L^2}+\|\mathcal{M}(H(f))\|_{L^2}.$$
Using explicit constants for $L^1$ norm of $\Psi$, the best constant for uncentered Hardy-Littlewood maximal operator, and unitarity of Hilbert transform, we have
$$\|H^*f\|_{L^2}\leq (1+\sqrt{2})\bigg(1+\frac{1+\log 2}{\pi}\bigg)\|Hf\|_{L^2}<\frac{15}{4}\|Hf\|_{L^2}.$$
For $L^{\infty}$ estimate we use that $\|\varphi_t\|_{L^1}=1$ and representation of $M^t$ as convolution operator. Again, by interpolation and the $L^{\infty}$ estimate of $H^*$, we finally get
$$\|H^*f\|_{L^p}<\bigg(\frac{15}{4}\bigg)^{\frac{2}{p}}\|Hf\|_{L^p},$$
for $p\geq 2.$
\end{remark}
\section{Appendix}\label{sec-4}
\subsection{Calculations of the integral (\ref{appendix})}\label{integral-37} In this section, we give some detailed calculations of the integral (\ref{appendix}).
By letting $\rho=2\pi r|x|$, we have
\begin{align*}
  \varphi_t(x)&=\frac{2\pi}{|x|^{\frac{d-2}{2}}}\int_{0}^{\infty}m(tr)J_{\frac{d-2}{2}}(2\pi r|x|)r^{\frac{d}{2}}\mathrm{d}r\\
  &=\frac{1}{(2\pi)^{\frac{d}{2}}|x|^d}\int_0^{\infty}m\left(\frac{t\rho}{2\pi|x|}\right)\mathrm{d}\left(\rho^{\frac{d}{2}}J_{\frac{d}{2}}(\rho)\right).
\end{align*}
It follows from (\ref{12-29-mx}) that
\begin{equation}
\label{1-9-eq2}
\frac{\mathrm{d}}{\mathrm{d}\rho}m\left(\frac{t\rho}{2\pi|x|}\right)=-\frac{2^{\frac{d}{2}}\Gamma(\frac{d+1}{2})}{\sqrt{\pi}}\left(\frac{t\rho}{|x|}\right)^{-\frac{d}{2}}J_{\frac{d}{2}}\left(\frac{t\rho}{|x|}\right)\frac{t}{|x|}.
\end{equation}
Therefore,
\begin{align*}
  \varphi_t(x)&=\frac{1}{(2\pi)^{\frac{d}{2}}|x|^d}m\left(\frac{t\rho}{2\pi|x|}\right)\left(\rho^{\frac{d}{2}}J_{\frac{d}{2}}(\rho)\right)\bigg|_{\rho=0}^{\infty} \\
   &\ \ \ +\frac{\Gamma(\frac{d+1}{2})t^{1-\frac{d}{2}}}{\pi^{\frac{d+1}{2}}|x|^{1+\frac{d}{2}}}\int_0^{\infty}J_{\frac{d}{2}}(\rho)J_{\frac{d}{2}}\left(\frac{t\rho}{|x|}\right)\mathrm{d}\rho.
\end{align*}

In what follows, we are going to show the first part of the above formula vanishes. Once this is done, then, by using (\ref{1-08-eq2}) and (\ref{4-26-eq1}), the desired result (\ref{1-9-eq1}) and (\ref{4-26-Im-2}) hold.
The proof is as follows:

(1) If $\rho=0$, then $m(0)=1$. It is easy to see that
$$\lim_{\rho\to0}m\left(\frac{t\rho}{2\pi|x|}\right)\left(\rho^{\frac{d}{2}}J_{\frac{d}{2}}(\rho)\right)=0,$$
since $d\geq2$ is an integer.

(2) If $\rho=\infty$, then the well-known asymptotic for oscillatory integral gives
\begin{equation}\label{9-16-2}\lim_{\rho\to\infty}J_{\frac{d}{2}}(\rho)=0.\end{equation}

Recall that (see (\ref{function-mx}))
$$m\left(\frac{t\rho}{2\pi|x|}\right)=\frac{2^{\frac{d}{2}}\Gamma(\frac{d+1}{2})}{\sqrt{\pi}}\int_{\frac{t}{|x|}\rho}^{\infty}r^{-\frac{d}{2}}J_{\frac{d}{2}}(r)\mathrm{d}r.$$
By using the following asymptotic estimate of Bessel functions (\cite[P. 195]{watson} and \cite[B. 8]{LG}):
$$J_{\frac{d}{2}}(r)=\sqrt{\frac{2}{\pi r}}\left(\cos\left(r-\frac{d+1}{4}\pi\right)+O\left(\frac{1}{r}\right)\right),$$
we have
$$m\left(\frac{t\rho}{2\pi|x|}\right) =\frac{2^{\frac{d+1}{2}}\Gamma(\frac{d+1}{2})}{\pi}\int_{\frac{t}{|x|}\rho}^{\infty}r^{-\frac{d+1}{2}}\left(\cos\left(r-\frac{d+1}{4}\pi\right)+O\left(\frac{1}{r}\right)\right)\mathrm{d}r.$$
The second part of the above integral can be estimated as follows:
$$\int_{\frac{t}{|x|}\rho}^{\infty}r^{-\frac{d+1}{2}}O\left(\frac{1}{r}\right)\mathrm{d}r=O\left(\int_{\frac{t}{|x|}\rho}^{\infty}r^{-\frac{d+3}{2}}\mathrm{d}r\right)=O\left(\rho^{-\frac{d+1}{2}}\right),\ \ \ \mbox{as}\ \ \rho\to \infty.$$
For the first part of the above integral, we can use integration by parts and obtain the following:
\begin{align*}
  \int_{\frac{t}{|x|}\rho}^{\infty}r^{-\frac{d+1}{2}}\cos\left(r-\frac{d+1}{4}\pi\right)\mathrm{d}r & =-\left(\frac{t}{|x|}\rho\right)^{-\frac{d+1}{2}}\sin\left(\frac{t}{|x|}\rho-\frac{d+1}{4}\pi\right)\\
   & \ \ \ +\frac{d+1}{2}\int_{\frac{t}{|x|}\rho}^{\infty}r^{-\frac{d+3}{2}}\sin\left(r-\frac{d+1}{4}\pi\right)\mathrm{d}r\\
   &=O\left(\rho^{-\frac{d+1}{2}}\right), \ \ \ \mbox{as}\ \ \ \rho\to\infty.
\end{align*}
Based on the above discussions, we have
\begin{equation}\label{9-16}m\left(\frac{t\rho}{2\pi|x|}\right)= O\left(\rho^{-\frac{1+d}{2}}\right),\ \ \ \mbox{as}\ \ \ \rho\to \infty.\end{equation}

Now, by using (\ref{9-16-2}) and (\ref{9-16}), one can easily obtain
$$\lim_{\rho\to\infty}m\left(\frac{t\rho}{2\pi|x|}\right)\rho^{\frac{d}{2}}J_{\frac{d}{2}}(\rho)=0,$$
which completes the proof. \qed

\subsection{Another computational proof of $\|\varphi_t\|_{L^1}=1$}\label{L1-norm-phi} 
Let $\varphi(x)$ be given by (\ref{function-phi}).
Then  $\varphi_t(x):=\varphi(x/t)/t^d$.
Since $\varphi$ is positive and radial, we see that $\|\varphi\|_{L^1}=\|\varphi_t\|_{L^1}$. This shows that the $L^1$ norm of $\varphi_t$ is independent of $t$.

Notice that $\|\varphi\|_{L^1}=\|\varphi_1+\varphi_2\|_{L^1}=\|\varphi_1\|_{L^1}+\|\varphi_2\|_{L^1}$,  if we define $\varphi_1=\varphi(x)\chi_{|x|>1}(x)$ and $\varphi_2=\varphi(x)\chi_{|x|<1}(x)$. 

We first calculate $\|\varphi_1\|_{L^1}$ as follows:
According to the definition of hypergeometric functions, we see that
\begin{equation}
\label{12-29-I1}\varphi_1(x)=\frac{1}{|x|^{1+d}\pi^{1+\frac{d}{2}}}\frac{\Gamma(\frac{1+d}{2})}{\Gamma(\frac{1}{2})}\sum_{n=0}^{\infty}\frac{\Gamma(n+\frac{1}{2})\Gamma(n+\frac{1+d}{2})}{\Gamma(n+1+\frac{d}{2})n!}\left(\frac{1}{|x|}\right)^{2n},\ \mbox{for}\ |x|>1.
\end{equation}
Moreover, it follows from \cite[Appendix A.3]{LG} that
$$\omega_{d-1}=\int_{\mathbb{S}^{d-1}}\mathrm{d}\theta=\frac{2\pi^{\frac{d}{2}}}{\Gamma(\frac{d}{2})}.$$
Then, since $\varphi_1(x)$ is radial and $\varphi_1(x)=0$ if $|x|\leq1$, by letting $x=r\theta$, we have
$$\|\varphi_1\|_{L^1}=\int_{\mathbb{R}^n}|\varphi_1(x)|\mathrm{d}x=\int_{1}^{\infty}\varphi_1(r)r^{d-1}\mathrm{d}r\int_{\mathbb{S}^{d-1}}\mathrm{d}\theta:=I_1\omega_{d-1},$$
where $I_1=\int_1^{\infty}\varphi_1(r)r^{d-1}\mathrm{d}r$.

Elementary calculations and using (\ref{12-29-I1}), we have
\begin{align*}
 I_1&=\frac{\Gamma(\frac{1+d}{2})}{\Gamma(\frac{1}{2})\pi^{1+\frac{d}{2}}}\sum_{n=0}^{\infty}\frac{\Gamma(n+\frac{1}{2})\Gamma(n+\frac{1+d}{2})}{\Gamma(n+1+\frac{d}{2})n!}\int_1^{\infty}\frac{1}{r^{2}}\left(\frac{1}{r}\right)^{2n}\mathrm{d}r \\
   &=\frac{\Gamma(\frac{1+d}{2})}{\Gamma(\frac{1}{2})\pi^{1+\frac{d}{2}}}\sum_{n=0}^{\infty}\frac{\Gamma(n+\frac{1}{2})\Gamma(n+\frac{1+d}{2})}{\Gamma(n+1+\frac{d}{2})n!(2n+1)}.
\end{align*}
Therefore,
\begin{align}\label{412-eq1}\nonumber
  \|\varphi_1\|_{L^1} &=\frac{2\pi^{\frac{d}{2}}}{\Gamma(\frac{d}{2})}\frac{\Gamma(\frac{1+d}{2})}{\Gamma(\frac{1}{2})\pi^{1+\frac{d}{2}}}\sum_{n=0}^{\infty}\frac{\Gamma(n+\frac{1}{2})\Gamma(n+\frac{1+d}{2})}{\Gamma(n+1+\frac{d}{2})n!(2n+1)} \\
   & =\frac{2\Gamma(\frac{1+d}{2})^2}{\pi\Gamma(\frac{d}{2})\Gamma(1+\frac{d}{2})}\ _3F_2\left[\frac{1}{2}, \frac{1}{2}, \frac{1+d}{2}; \frac{3}{2}, 1+\frac{d}{2}; 1\right].
\end{align}

Similarly, since
$$\varphi_2(x)=\frac{\Gamma(\frac{1+d}{2})}{\pi^{\frac{d+3}{2}}}\sum_{n=0}^{\infty}\frac{\Gamma(n+\frac{1}{2})\Gamma(n+\frac{1+d}{2})}{\Gamma(n+1+\frac{d}{2})n!}|x|^{2n}, \ \mbox{for}\ |x|<1,$$
is radial, and $\varphi_2=0$ if $|x|\geq1$, we have
\begin{align}\label{4-29-eq1}\nonumber
  \|\varphi_2\|_{L^1} &=\frac{2\Gamma(\frac{1+d}{2})}{\pi^{\frac{3}{2}}\Gamma(\frac{d}{2})}\sum_{n=0}^{\infty}\frac{\Gamma(n+\frac{1}{2})\Gamma(n+\frac{1+d}{2})}{\Gamma(n+1+\frac{d}{2})n!(2n+d)} \\
   & =\frac{2\Gamma(\frac{1+d}{2})^2}{\pi\Gamma(\frac{d}{2})\Gamma(1+\frac{d}{2})}\frac{1}{d}\ _3F_2\left[\frac{1}{2}, \frac{d}{2}, \frac{1+d}{2}; 1+\frac{d}{2}, 1+\frac{d}{2}; 1\right].
\end{align}

By using \cite[Corollary 2.4.5]{Landrews}, one has
\begin{align*}
  \ _3F_2\left[\frac{1}{2}, \frac{d}{2}, \frac{1+d}{2}; 1+\frac{d}{2}, 1+\frac{d}{2}; 1\right] & =\frac{\pi\Gamma(\frac{d}{2}+1)^2}{\Gamma(\frac{d+1}{2})^2}- \\
   & \frac{2\Gamma(\frac{d}{2}+1)^2}{\Gamma(\frac{d}{2})\Gamma(\frac{1}{2})\Gamma(\frac{d+3}{2})}\ _3F_2\left[1,1,\frac{1+d}{2};\frac{3}{2}, \frac{d+3}{2}; 1\right]
\end{align*}
and (by the identity \cite[(16.4.11)]{HB-10})
$$\ _3F_2\left[\frac{1}{2}, \frac{1}{2}, \frac{1+d}{2}; \frac{3}{2}, 1+\frac{d}{2}; 1\right]=\frac{\Gamma(\frac{d}{2}+1)}{\Gamma(\frac{1}{2})\Gamma(\frac{d+3}{2})}\ _3F_2\left[1,1,\frac{1+d}{2}; \frac{3}{2}, \frac{d+3}{2}; 1\right].$$
Now, putting these two expressions into (\ref{412-eq1}) and (\ref{4-29-eq1}), we get the desired result
$\|\varphi_t\|_{L^1}=\|\varphi\|_{L^1}=\|\varphi_1\|_{L^1}+\|\varphi_2\|_{L^1}=1.$
\qed


\vspace*{5mm}
\noindent {\bf Acknowledgments}.
We would like to thank Professor Blazej Wr\'obel for his helpful comments and suggestions on this paper.
We would also thank the anonymous referee for his/her helpful comments that have significant impact on this paper.

\vspace*{5mm}

\noindent {\bf Funding}.
The research of the first author was supported by NSFC (Grant No. 11925107, 12226334), the second author was supported by NSFs of Serbia (MPNTR 174017), and the third author was supported by NSFC (Grant No. 12271189, 11971182), NSF of Fujian Province (Grant No. 2021J01304, 2023J01127).

\vspace*{5mm}
\noindent {\bf Conflict of Interests}.
The authors declare that there is no conflict of interests regarding the publication of this paper.

\vspace*{5mm}
\noindent {\bf Data Availability Statement}.
The authors declare that this research is purely theoretical and does not
associate with any data.

\end{document}